\documentclass{article}

\usepackage{hyperref}
\usepackage{cite}
\usepackage{enumitem}
\usepackage[ruled,longend]{algorithm2e}
\usepackage{todonotes}

\usepackage{fancyhdr}
\usepackage{caption}
\usepackage{subcaption}
\usepackage{float}
\usepackage{amsmath,amssymb,amsthm,graphicx,url}

\theoremstyle{plain}
\newtheorem{lemma}{Lemma}[section]
\newtheorem{proposition}[lemma]{\textbf{Proposition}}

\newtheorem{theorem}[lemma]{\textbf{Theorem}}

\theoremstyle{definition}
\newtheorem{definition}[lemma]{\textbf{Definition}}
\newtheorem{example}[lemma]{\textbf{Example}}

\newtheorem{remark}[lemma]{Remark}

\title{A graphical algorithm for the integration of monomials 
in the Chow ring of the moduli space of stable marked curves of genus zero}

\author{Jiayue Qi \thanks{Doctoral Program ``Computational Mathematics'' W1214, Johannes Kepler 
University Linz.}}
\date{}
\begin{document}
\maketitle

\abstract{The Chow group of zero cycles in the moduli space of stable pointed curves of genus zero
is isomorphic to the integer additive group. Let $M$ be monomial in this Chow group. If no two factors
of $M$ fulfill a particular quadratic relation, then the monomial can be represented equivalently by a
specific tree; otherwise, $M$ is mapped to zero under the stated isomorphism.
Starting from this tree representation, we introduce a graphical algorithm for computing the 
corresponding integer for $M$ under the aforementioned isomorphism. The algorithm is linear with respect to the 
size of the tree.}
\section{Introduction}

The moduli space of stable $n$-pointed curves of genus zero, denoted by $\mathcal{M}_n$ in our paper, 
is a renowned object in modern intersection theory; for example,
it is the base for the definition of Gromov-Witten invariants~\cite{GW}. It is a smooth irreducible projective 
compactification of stable $n$-pointed genus-zero curves, which was originally introduced 
by Knudsen and Mumford in their series of papers \cite{km1}, \cite{km2} and \cite{km3}. 

Chow rings are essential in intersection theory, to indicate 
the intersection numbers of subvarieties. Let $X$ be a projective variety 
of dimension $k$. The Chow ring $A^*(X)$ of $X$ is a graded ring, and is the direct 
sum of $k+1$ groups, each of them is composed of cycles (formal sums of subvarieties)
of a fixed dimension. Conventionally, the group constituted by cycles of codimension
$r$ is called the Chow group of codimension $r$ and denoted by $A^r(X)$, and $A^r(X)=\{0\}$ when $r>k$. 
Particularly, $A^k(X)$ is the Chow group of cycles of dimension zero and is isomorphic to the integer 
additive group. 

In this paper, we work in the Chow ring of the moduli space $\mathcal{M}_n$ of stable $n$-pointed curves of genus zero
and we denote it by $A^*(n)$. Since $\mathcal{M}_n$ is of dimension $n-3$, we know that 
$A^*(n)=\bigoplus^{n-3}_{r=0}{A^r(n)}$, where $A^r(n)$ denotes the Chow group of codimension $r$.
We have $A^r(n)=\{0\}$ when $r>n-3$, and $A^{n-3}(n)\cong \mathbb{Z}$ --- we use the notation $\int$
to denote this isomorphism, following the convention. The integer under this isomorphism is 
called the {\em integral value} or {\em value} of the given element in $A^{n-3}(n)$.  
A set of generators of the group $A^1(n)$ was given in Keel's paper \cite{keel}, where each generator 
is indexed by a bi-partition of $\{1,\cdots,n\}$ and this set is also the generating set
for the whole ring. 

We will introduce an algorithm for computing the integral value of a product of the Keel generators, that is,
a monomial in the generators. The problem originally showed up as a sub-problem for counting the realization
of Laman graphs (minimally-rigid graphs) on a sphere \cite{laman}, when we wanted to improve the algorithm
given in \cite{laman}. With the help of this algorithm, we invented another algorithm for the same goal. However, by efficiency 
it does not seem faster than the one provided in \cite{laman}. But we see this problem fundamental, 
standing on its own and we find the algorithm elegant and concise, also
may be helpful for other similar or even further-away problems. Therefore, we formulate
it on its own.
We consider the situation when this monomial is of degree $n-3$, otherwise
we define its value to be zero. The input is $n-3$ such generators, hence is linear in $n$ and the output 
is an integer. Our algorithm is polynomial in $n$. 

A quadratic relation between the generators was introduced in Keel's paper \cite{keel}. This relation is called 
{\em Keel's quadratic relation} in our paper. With the help of this 
relation, we know that if any two factors of fulfills the relation, then the whole monomial has value zero. 
This observation naturally formed the first step of our algorithm. We check if any of the two factors of the input
monomial fulfills this relation: if yes, return zero; otherwise, we consider an equivalent characterization 
of the given monomial, in a specific tree --- {\em loaded tree}. This first part is polynomial in $n$ in the worst case.
In the second step, we transfer this tree 
via three steps to a forest. Next, we compute the integral value of the given monomial directly from the obtained forest.
The second part is linear with respect to $n$.

The tree characterization of the given monomial is inspired by 
\cite[Section~2.2.]{QS:representation} and was first introduced in \cite{ACM_communication}. Our paper can be considered
as a proper extension for \cite{ACM_communication}. We will describe the same algorithm, but with all the
detailed proofs provided. Note that the same problem was considered in \cite{vertex_splitting}, where an 
equivalent characterization for the algebraic reductions (in the ambient ring) for the input monomial is given, as some operations
on the tree representation. In fact, we were motivated by that characterization, then
we started to try-out many examples using the algorithm provided in \cite{vertex_splitting}. Eventually and excitingly, 
we discovered our algorithm which as an algorithm, is much more efficient and neat, comparing to the one given
in \cite{vertex_splitting}. 
Later on, we managed to prove the correctness of our algorithm with the help of some more advanced algebraic geometry tools.

\section{Preliminaries}
Since the main problem this paper focus on is exactly the same with the papers \cite{ACM_communication} and \cite{vertex_splitting},
the preliminaries will be very much similar to them. However, for completeness, we introduce everything from scratch.

Let $n\in \mathbb{N}$, $n\geq 3$, define $N:=\{1,\dotsc,n\}$. Denote by
$\mathcal{M}_n$ the moduli 
space of stable $n$-pointed curves of genus zero. 
A bi-partition $\{I,J\}$ of $N$ where the cardinalities of $I$ and $J$ are both at least
$2$ is called a {\bf cut}; we call $I$ and $J$ {\bf parts} of this cut.
There is a hypersurface $D_{I,J}\in \mathcal{M}_n$ for each cut $\{I,J\}$ and its class in the Chow ring 
is denoted by $\delta_{I,J}$. Note that $D_{I,J}=D_{J,I}$, and as well $\delta_{I,J}=\delta_{J,I}$.    
This Chow ring is a graded ring and we denote it by $A^*(n)$. Then we have 
$A^*(n)=\bigoplus^{n-3}_{r=0}{A^r(n)}$.
These homogeneous components are defined as the Chow groups (of $\mathcal{M}_n$); $A^r(n)$ is the {\bf Chow group of codimension $r$}. 
It is known that $A^r(n)=\{0\}$ for $r > n-3$ and
 $A^{n-3}(n)\cong \mathbb{Z}$. We denote this isomorphism by
$\int: A^{n-3}(n) \to \mathbb{Z}$. We can extend this map so that it is defined on the whole ring by setting the value of 
all other elements to be zero:
$\int:A^*(n)\to \mathbb{Z}$, $\int(M)=0$ if $M\notin A^{n-3}(n)$; it is then a group homomorphism between the Chow ring and 
integer additive group, we call it the {\em integral map}.

The set $\{\delta_{I,J}\, \mid \, \{I,J\} \text{ is a cut}\}$ generates the group $A^1(n)$, and also the 
whole ring $A^*(n)$. Then we can view
$\prod_{i=1}^{n-3}{\delta_{I_i,J_i}}$ as an element 
in $A^{n-3}(n)$, since $A^*(n)$ is a graded ring. We define the {\bf value} of $M:= \prod_{i=1}^{n-3}{\delta_{I_i,J_i}}$ to be  
$\int(\prod_{i=1}^{n-3}{\delta_{I_i,J_i}})$. 
The problem we deal with in this paper
is to {\bf compute the value of a given monomial $M= \prod_{i=1}^{n-3}{\delta_{I_i,J_i}}$.}

Keel introduced a quadratic relation between the generators of $A^*(n)$ in \cite{keel}; 
we call it {\em Keel's quadratic relation}. 
We say that two generators
 $\delta_{I_1,J_1}, \delta_{I_2,J_2}$ fulfill {\bf Keel's quadratic relation} (\cite[Section 4, Theorem~1.(3)]{keel}) if the following 
 four conditions hold:
 $I_1\cap I_2\neq \emptyset$;
$I_1\cap J_2\neq \emptyset$;
 $J_1\cap I_2\neq \emptyset$;
$J_1\cap J_2\neq \emptyset$.
And in this case, we have $\delta_{I_1,J_1}\cdot \delta_{I_2,J_2}=0$. For example, when $n=5$,
$\delta_{12,345}\cdot \delta_{14,235}=0$ since these two factors fulfill the Keel's quadratic relation. 
Note that we use abbreviated notations for the index
of the generators, for instance $\delta_{12,345}$ represents $\delta_{\{1,2\},\{3,4,5\}}$. We
will use this abbreviation also in the later context.
Motivated by this quadratic relation, we realize that when ever two factors of the given monomial
fulfill this relation, the value of the monomial is zero. 
Hence the first step of our algorithm 
is to check whether there are two factors fulfilling this relation. There are in total $n-3$ input 
factors, so we need to check ${n-3 \choose 2}$ many pairs in the worst case. The checking for each
pair of generators is linear in $n$. Hence the algorithm in this step is polynomial in the worst case.

We call those monomials of which no two factors fulfill the Keel's quadratic relation {\bf tree monomials},
since there exists a one-to-one correspondence between these monomials and a specific type of trees which 
we call {\em loaded trees} (see Definition \ref{def:loaded_tree}).
Then, how should we compute the value of a tree monomial (in $A^{n-3}(n)$)? The following theorem 
indicates the first thing to check, when we have a tree monomial at hand.
\begin{theorem}\cite[Theorem 3.4.]{vertex_splitting}\label{thm:clever}
 If all factors are distinct in the tree monomial $M:=\Pi^{n-3}_{i=1}{\delta_{I_i,J_i}}$, where 
 $I_i\cup J_i=N$ for all $1\leq i\leq n-3$.
 Then $\int(M)=1$ and we call this type of
 tree monomials {\bf clever monomials}.
\end{theorem}
Actually, we can combine the first two steps as one step, serving as the first part of our algorithm. This is because
going through all the pairs of generators once is sufficient: we can check whether the pair fulfills the 
Keel's quadratic relation or not and whether the generators in the pair are distinct at the same time. Hence the first part of our algorithm 
is polynomial with respect to $n$.

\subsection{Loaded trees}

In this section, we introduce the one-to-one correspondence between tree monomials and loaded trees.

\begin{figure}
\centering
\begin{subfigure}{.4\textwidth}
  \centering
  \includegraphics[width=0.8\linewidth]{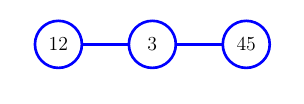}
  \caption{}
\end{subfigure}%
\begin{subfigure}{.4\textwidth}
  \centering
  \includegraphics[width=0.8\linewidth]{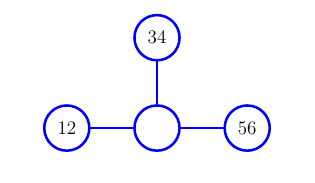}
  \caption{}
\end{subfigure}
\caption{On the left is a loaded tree with $5$ labels and $2$ edges; its monomial is $\delta_{12, 345}\cdot \delta_{123,45}$.
On the right is a loaded tree with $6$ labels and $3$ edges; its monomial is 
$\delta_{12,3456}\cdot \delta_{34,1256}\cdot \delta^2_{56,1234}$. Note that we use an abbreviated notation
for the labeling set of vertices, for instance $12$ refers to $\{1,2\}$, for a neater picture.}
\label{fig:loaded_trees}
\end{figure}

\begin{definition}[\cite{ACM_communication}, Definition 0.1.]\label{def:loaded_tree}
 A {\bf loaded tree with $n$ labels and $k$ edges} is a tree $T=(V,E)$ together with a 
 labeling function $h: V\to 2^N$ and an edge multiplicity function $m: E \to \mathbb{N}^+$ such that the
 following three conditions hold:
 \begin{enumerate}
  \item $\{h(v)\}_{v\in V, h(v)\neq \emptyset}$ form a partition of $N$; elements in $N$ are called the {\em labels}
  of $T$.
 \item $\sum_{e\in E}{m(e)}=k$.
 \item For every $v\in V$, $\deg(v)+|h(v)|\geq 3$, note that here multiple edges are only counted once for the degree 
 of its incident vertices.
\end{enumerate}
\end{definition}
We define the {\em monomial of a loaded tree} as follows. Removing any edge $e$ gives us two components of the tree, 
then the two sets of labels in each components respectively gives us a cut $\{I,J\}$ of $N$; $\delta_{I,J}$ is the 
corresponding factor of the edge $e$. The product of the corresponding factors of all edges is defined to be the monomial
of the given loaded tree. We see that a loaded tree uniquely determines its monomial. 
We see two examples of loaded trees and their monomials in Figure \ref{fig:loaded_trees}. 
Note that in the example, we use an abbreviated notation for 
the labeling set of vertices shown on the picture. We will keep using it in the later context, for neater pictures.

The following theorem tells us the existence of a one-to-one correspondence between tree monomials and loaded trees. 
 \begin{theorem}\cite[Theorem 2.2.]{vertex_splitting}\label{thm:correspondence}
  There is a one-to-one correspondence between tree monomials 
 $M=\prod^k_{i=1}{\delta_{I_i,J_i}}$ and loaded trees with $n$ labels and $k$ edges, where $I_i\cup J_i=N$ for all 
 $1\leq i\leq k$. 
 \end{theorem}
\begin{remark}
 Note that the original idea of this correspondence came from Section 2.2 of paper \cite{QS:representation}.
 However, it is better explained and stated in \cite{vertex_splitting}, in the sense of monomials whose factors are 
 indexed by cuts of $N$, i.e., generators of $A^*(n)$. Usually we denote by $M_T$ the monomial of loaded tree $T$, 
 and by $T_M$ the loaded tree of monomial $M$. We call the corresponding loaded tree for clever monomials {\bf clever trees}.
 \end{remark}

It is trivial to obtain the monomial of a loaded tree, while the other direction not. The algorithm for this direction
is described in \cite{ACM_communication}. Note that it specifies the ambient group of the input monomial, but
the same algorithm also works for a monomial of any other degree. The idea of this algorithm comes from Section 2.2 of
the paper \cite{QS:representation}. We will not go into details of this algorithm in the paper.
For completeness, we also illustrate this algortihm here, see 
Algorithm \ref{alg:monomial_to_tree}.

 \begin{algorithm}
\thispagestyle{empty}
\caption{monomial to tree (Algorithm 1 in \cite{ACM_communication})}
\label{alg:monomial_to_tree}
\SetKwInOut{Input}{input}
\SetKwInOut{Output}{output}

\Input{a tree monomial $M$ in $A^{n-3}(n)$}
\Output{a loaded tree with $n$ labels and $n-3$ edges}
$C \gets$ collection of any cut that corresponds to some factor of $M$\;
$P \gets$ collection of all the parts of cuts in $C$\;
$c \gets$ any element $c=\{I,J\}\in C$\;

\For{each element $p\in P\setminus \{I,J\}$}
   {\If{$p\subset I$ or $p\subset J$}
       {$c:=c\cup \{p\}$}}
$H \gets$ the Hasse diagram of elements in $c$ with respect to set containment order\;
Consider $H$ as a graph $(V,E)$\;
\For{each vertex $V$ of $H$}
    {Define labelling set $h(V)$ as its corresponding element in $c$\;
    Update the labelling set: $h(V):=h(V)\setminus h(V_1)$ if $V_1$ is less than $V$ in $H$}

$E:=E\cup \{\{I,J\}\}$\;
Attach this labelling function $h$ to $H$\;
Set the multiplicity function value $m(e)$ for each edge $e$ 
as the power of its corresponding factor in $M$\;

\Return{$H=(V,E,h,m)$}         
  
\end{algorithm}

Let us see an example on constructing the corresponding loaded tree of a given monomial, using Algorithm \ref{alg:monomial_to_tree},
so as to have an intuitive comprehension.
\begin{example}\label{eg:monomial_to_tree}
 Given a tree monomial $\delta^3_{123,456789}\cdot \delta_{12345,6789}\cdot \delta_{1234589,67}\cdot \delta_{1234567,89}$. 
 Obviously we have the set of labels $N:=\{1,2,3,4,5,6,7,8,9\}$.
 We collect the parts in set 
 $P:= \{\{1,2,3\},\{4,5,6,7,8,9\},\{1,2,3,4,5\},\{6,$ 
 $7,8,9\},\{1,2,3,4,5,8,9\},\{6,7\},\{1,2,3,4,5,6,7\},\{8,9\}\}$
 and we pick any cut $c=\{\{1,2,3,4,5\},\{6,7,8,9\}\}$
 from the set of cuts. After collecting all parts which are either contained in $\{1,2,3,4,5\}$ or $\{6,7,8,9\}$, 
 we obtain $c=\{\{1,2,3,4,5\},\{6,7,8,9\},\{1,2,3\},\{6,7\},\{8,9\}\}$, then we construct the corresponding Hasse diagram for $c$, 
 see Figure \ref{fig:hasse}. The output loaded tree $T_M$ is shown in Figure \ref{fig:from_hasse}. 
 It is easy to see that if we go back from the tree constructing monomial, we again obtain $M$.
 
  \begin{figure}[H]
\centering
\begin{subfigure}{.4\textwidth}
  \centering
  \includegraphics[width=0.8\linewidth]{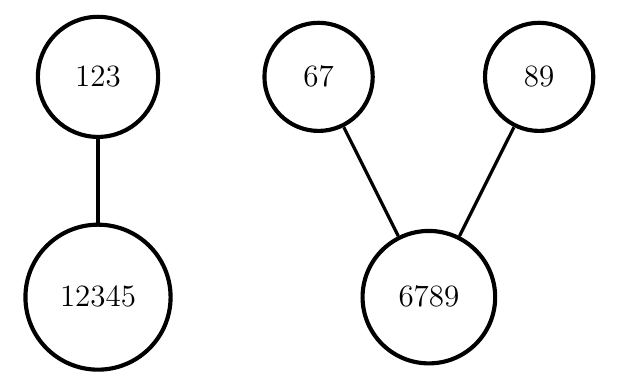}
  \caption{}
  \label{fig:hasse}
\end{subfigure}%
\begin{subfigure}{.4\textwidth}
  \centering
  \includegraphics[width=0.9\linewidth]{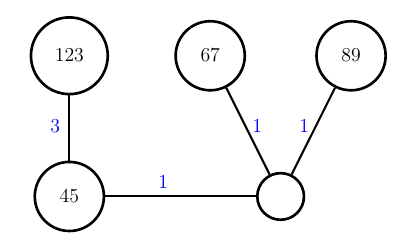}
  \caption{}
  \label{fig:from_hasse}
\end{subfigure}
\caption{On the left is the Hasse diagram of set $\{\{1,2,3,4,5\},\{6,7,8,9\},$ 
  $\{1,2,3\},\{6,7\},\{8,9\}\}$ with respect to set containment order. On the right is 
  the corresponding loaded tree of monomial 
$\delta^3_{123,456789}\cdot \delta_{12345,6789}\cdot \delta_{1234589,67}\cdot \delta_{1234567,89}$. 
Multiplicity function values are written in blue.}
\end{figure}
 \end{example}
 
 Note that if a loaded tree has no edges, then its monomial has no factors; we
call such monomial an {\bf empty monomial}. We extend this one-to-one correspondence a little by including 
a single tree and a single monomial: the loaded tree with $3$ labels and no edges corresponds to
the empty monomial; the loaded tree with $k$ labels and no edges has no corresponding monomial if $k\neq 3$.
With this extension proposed, we can now define the {\bf value of a loaded tree} as the value of its corresponding 
monomial. Hence our goal can be expressed in other words now: {\bf compute the value of a loaded tree with $n$ labels and 
$n-3$ edges, where $n\geq 3$.} We call such loaded trees {\bf proper loaded trees}.
Later we will see, the extension above is done so as to guarantee that the loaded tree has the same value 
with its monomial, which also ought to hold by definition. 

This tree representation is the foundation for our algorithm, and constitute the second part of our algorithm ---
transferring the tree monomial to its corresponding loaded tree; this step is at most polynomial to $n$.
However, all the main contents so far are already mentioned
in \cite{vertex_splitting}. In the next section, we introduce our graphical algorithm on computing the value of a loaded tree,
i.e., the third part of our algorithm, which chiefly reflect the efficiency and conciseness of our algorithm.

\section{The graphical algorithm}
In this section, we illustrate a graphical algorithm, computing the value of a proper loaded tree. We already know that 
the value of a clever tree is one --- this will just be a special case for our general algorithm.
We postpone the correctness proof of the algorithm to later sections. Note that the algorithm is the same as 
described in \cite{ACM_communication}, only some terms (names of the trees) are modified. 

Here we give the sketch of the third part of our algorithm:
\begin{enumerate}
\item {\bf Input:} a loaded tree with $n$ labels and $n-3$ edges, i.e. a proper loaded tree with $n$ labels.
\item {\bf Output:} the value of the input loaded tree (which is an integer).
\item Transfer the loaded tree to a \textbf{weighted tree}.
\item Calculate the \textbf{sign of the tree value}.
\item Construct a \textbf{redundancy tree} from the weighted tree.
\item Construct a \textbf{redundancy forest} from the redundancy tree.
\item Apply a recursive algorithm to each tree in the redundancy forest, obtaining 
 the absolute tree value.
\item Product of the sign and absolute value gives us the tree value.
\end{enumerate}
In the sequel, we will explain the algorithm step by step, with a running example.

A {\bf weighted tree} is a tuple $(V,E,w)$ where $(V,E)$ is a tree and $w:V\cup E\to \mathbb{N}$
is the weight function assigning to each edge and each vertex a non-negative integer weight. 
Let $LT=(V,E,h,m)$ be a loaded tree, define a function $w:V\cup E\to \mathbb{N}$ by $w(v):=\deg(v)+|h(v)|-3$
for all $v\in V$ and $w(e):=m(e)-1$ for all $e\in E$. From the third item of Definition \ref{def:loaded_tree},
we see that $w(v)\geq 0$ for all $v\in V$ and naturally multiplicity of any edge is in $\mathbb{N}^+$.
Hence, $(V,E,w)$ is a weighted tree; we call it {\em the weighted tree of $LT$}. 
It is not hard to verify that $\sum_{v\in V}{w(v)}=\sum_{e\in E}{w(e)}$ holds for a weighted tree
if it is of some proper loaded tree. This
identity is called the {\em weight identity} \cite[Section 2.]{vertex_splitting}.
From Remark 3.8. and the reduction chain algorithm (Section 5.2) of \cite{vertex_splitting},
we know that if the value of the given loaded tree is non-zero, then {\bf the sign of the tree value
is $-1$ to the power of the edge weight sum (or equivalently, the vertex weight sum)}. 
This is because each recursion step contributes a negative sign to the value, and from the linear reduction
(\cite[Section 3]{vertex_splitting}) we know that each recursive step reduces the edge weight sum by one. 

We say that the tuple $(V,E,w)$ is a {\bf redundancy tree} if $(V,E)$ is a tree and $w:V\to \mathbb{N}$
is a function defined on the vertex set. In the last step of the algorithm, we obtained a 
weighted tree. Then, we replace each edge by two edges with a vertex in the middle ---
inheriting the weight of the replaced edge --- connecting them. We see that in this way, we actually
get a redundancy tree. We call the so-gained tree the {\em redundancy tree of the given loaded tree (or, of the given weighted tree)}.
A {\bf redundancy forest} is defined to be a forest in which each tree is a redundancy tree.
From the redundancy tree we obtained in the last step, we will obtain a redundancy forest via 
deleting all vertices of zero weight-value and their incident edges; this so-obtained forest is then 
called {\em the redundancy forest of the given loaded tree / weighted tree / redundancy tree}.

After we obtain the redundancy forest of the given loaded tree, we apply a recursive formula to each redundancy tree in the forest,
so as to obtain the absolute value of the loaded tree. Let $RF$ be the redundancy forest of loaded 
tree $LT$, define the {\bf value of $RF$} (denoted by $\int(RF)$) as the product of the values of all the trees in the forest.
Define the {\bf value of a redundancy tree $RT=(V,E,w)$} recursively as follows.
Pick any leaf $l\in V$, compare the 
weight of $l$ and that of its unique parent $l_1$: if $w(l)>w(l_1)$, then return $0$; otherwise, 
$\int(RT):={w(l_1)\choose w(l)}\cdot \int(RT_1)$, where $RT_1=(V_1,E_1,w_1)$ is the redundancy tree defined as 
follows. Deleting leaf vertex $l$ and its incident edge from $RT$, then replace the weight of $l_1$ by $w(l_1)-w(l)$.
Formally speaking, we have $V_1=V\setminus \{l\}$, $E_1=E\setminus \{l,l_1\}$, $w_1(l_1)=w(l_1)-w(l)$ and $w_1(v)=w(v)$
for all $v\in V_1\setminus \{l_1\}$. When $RT$ is a degree-zero vertex, $\int(RT)=0$ if it has non-zero weight and 
$\int(RT)=1$ otherwise. If $RT$ is a null graph ---  the graph that contains no vertices or edges --- then 
$\int(RT)=1$. 
Then we say that the value of loaded tree $LT$ is the value of the redundancy forest $RF$.

Let us see an example, on how to obtain the value of a given loaded tree.

\begin{example}\cite[Example 0.4.]{ACM_communication}
Figure \ref{fig:loaded_tree} depicts $LT$ --- a loaded tree with $14$ labels and $11$ edges,
while Figure \ref{fig:weighted_tree} shows the weighted tree $WT$ of $LT$.
We obtain that the edge weight sum of $WT$ is $2+4+0+1=7$, while its vertex weight sum
$1+1+4+0+1=7$. Then we obtain that the sign of $\int(LT)$ is $(-1)^7=-1$.
Figure \ref{fig:redundancy_tree} shows the redundancy tree $RT$ of $LT$ (or of $WT$), 
and Figure \ref{fig:redundancy_forest} describes the corresponding redundancy forest $RF$.
Then apply the recursive formula on $RF$, we obtain 
$$\int(RF)=
 [{1 \choose 1} \times 1] \times [{4 \choose 1}\times {4 \choose 3}\times {2 \choose 1} \times {1 \choose 1} \times 1]=32.$$ 
 Product of $32$ and $-1$ tells us that the value of the loaded tree $LT$ shown in Figure \ref{fig:loaded_tree} is $-32$.

  \begin{figure}
\centering
\begin{subfigure}{0.5\textwidth}
  \centering
  \includegraphics[width=0.8\linewidth]{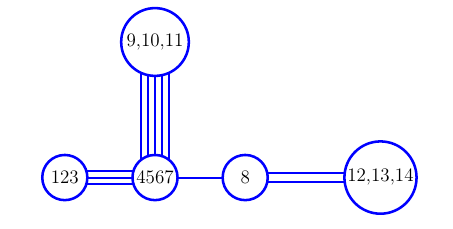}
  \caption{}
  \label{fig:loaded_tree}
\end{subfigure}%
\begin{subfigure}{0.5\textwidth}
  \centering
  \includegraphics[width=0.9\linewidth]{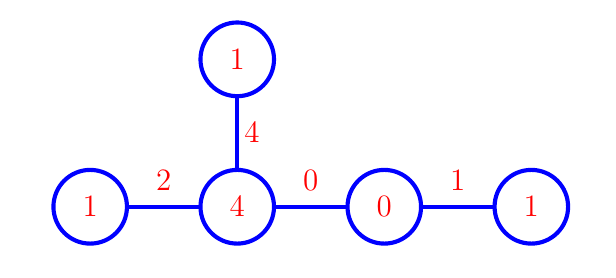}
  \caption{}
  \label{fig:weighted_tree}
\end{subfigure}
\caption{On the left is a loaded tree $LT$ with $14$ labels and $11$ edges. Labels are tagged in black.
On the right is the weighted tree $WT$ of the loaded tree $LT$ described in Figure \ref{fig:loaded_tree}. 
Weight function values are marked in red. }
\end{figure}

  \begin{figure}
\centering
\begin{subfigure}{0.5\textwidth}
  \centering
  \includegraphics[width=0.8\linewidth]{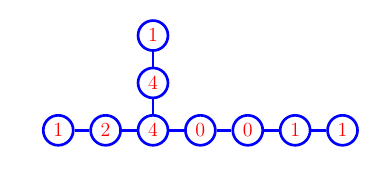}
  \caption{}
  \label{fig:redundancy_tree}
\end{subfigure}%
\begin{subfigure}{0.5\textwidth}
  \centering
  \includegraphics[width=0.9\linewidth]{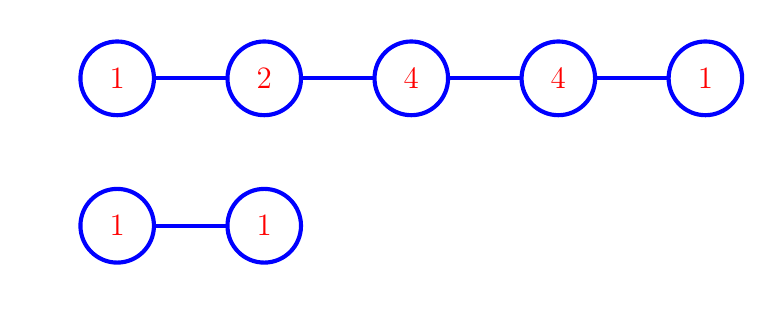}
  \caption{}
  \label{fig:redundancy_forest}
\end{subfigure}
\caption{On the left is the redundancy tree $RT$ of $LT$ described in Figure \ref{fig:loaded_tree}
or, of the weighted tree $WT$ depicted in Figure \ref{fig:weighted_tree}. Weight values are marked in red.
On the right is the redundancy forest obtained from $RT$ shown on the left.}
\end{figure}

\end{example}

Consider the whole procedure, from inputting a proper loaded tree, to finally obtaining the tree value. Termination is trivial since 
the input tree is finite, so does its redundancy forest.
In the recursion formula calculus, each step strictly reduces the size of the redundancy forest. 
Because of the following identity of binomial coefficients, we know that the recursive formula gives 
us the same value, no matter in which sequence do we consider and delete the leaf vertices of a redundancy tree.
$${c\choose a}\cdot {c-a \choose b} \equiv {c\choose b}\cdot {c-b \choose a}.$$
Therefore, the above process is indeed an algorithm. We call it the {\bf forest algorithm}. It is not hard to see that
the complexity of the forest algorithm is linear with respect to the number of vertices of the 
input loaded tree.
If the input loaded tree is a clever tree, then any vertex of its redundancy tree has value one.
Therefore, its redundancy forest is a null graph and hence the input loaded tree has value one. 
This indicates that the value of clever trees can also be handled by the forest algorithm, as a special case.

Based on the forest algorithm, let us consider again the extension of the one-to-one correspondence between loaded trees and tree monomials:
a loaded tree with a single vertex and $3$ labels has the null graph as its redundancy forest, hence has value one;
a loaded tree with a single vertex and $k$ ($k\neq 3$) labels has the single vertex with non-zero weight as its 
redundancy forest, hence has value zero. 
And the correctness on the base cases directly come from the definition of $\int$.
Given the correctness of the forest algorithm, our previous extension stands defensible.
In the next section, we prove the correctness of the forest algorithm.

\section{Correctness}
This section contains three sub-sections. First we introduce three kinds of edge-cutting operation on trees, 
and describe the main theorem using algebraic language. 
Second, we prove the correctness of the forest algorithm using an algebraic theorem posted in the first part.
Third, we use pure algebraic tools to prove the main theorem.

\subsection{Three types of edge-cutting}
In this section, we introduce three different types of edge-cutting, and using these graphical operations,
express the main theorem in algebraic language.
In order not to interrupt the story later on, we need to introduce the concept of {\em star-cut} first.
 \begin{definition}[star]
 Let $T=(V,E)$ be a tree with $|V|\geq 3$. If there exists $v\in V$ all other vertices is a neighbor of $v$, 
 then we call this tree a {\bf star}.
\end{definition}
Given a tree $T=(V,E)$. We can pick any edge $e=\{u,v\}\in E$. We call the process of cutting off edge $e$, attaching a new vertex
$u_1$ to $u$, $v_1$ to $v$ via a new edge $e_1$ and $e_2$, respectively
an {\bf edge-cut} of $T$; after 
this process, we obtain two new trees $T_1$, $T_2$ from $T$. If we obtain a star after applying edge cut on some edge of $T$, we 
call this edge cut a {\bf star-cut}.

\begin{proposition}
 Star cut exists for any tree with no less than three vertices.
\end{proposition}
\begin{proof}
 Let $T$ be an arbitrary tree with no less than three vertices. Define $L$ to be the set of leaves of $T$. 
 Then define an equivalence relation on $L$ by
 $v_1\sim v_2$ iff $N(v_1)=N(v_2)$ for any two leaves $v_1,v_2$ of $T$, where $N(v)$ refers to the set of neighbors of vertex $v$. 
 It is not hard to see that there is a 1-1 correspondence 
 between support (non-leaf) vertex set of $T$ and set of equivalence classes we define above. 
 
If $T$ is a star, then the proposition holds
 since we can apply edge-cut to any edge of $T$ and we will get a star after it. 
 Otherwise we delete all leaves of $T$, obviously we obtain a nontrivial
 tree $T_1$ --- here non-trivial means that it is not a single vertex. So it
 must have an vertex $u$ with degree $1$, w.l.o.g., assume $\{u,u'\}\in E(T_1)$, then we also have $\{u,u'\}\in E(T)$. 
 Obviously, $u$ is a support vertex of $T$. If we apply edge-cut to edge $\{u,u'\}$ in $T$, we get a star centered at 
 vertex~$u$, it is a star-cut.
 \end{proof}
 
We continue by recalling the {\bf single-edge cutting} operation defined in \cite{vertex_splitting}.
Let $T$ be a loaded tree and $e=\{u,v\}$ be any edge of $T$ with multiplicity $r$ and corresponding factor
$\delta_{I_1,I_2}$. When $r=1$, we construct two other loaded trees $T_1$ 
and $T_2$ by cutting off edge $e$ and add one more label $x$ to $u$ and $v$, respectively.
Imagine that we cut off edge $e$, we would obtain two components: all the labels in one component form the set $I_1$; we call 
this component {\em Component-$I_1$} and the other component {\em Component-$I_2$}.
Assume that $u$ is in Component-$I_1$ and $v$ is in Component-$I_2$.
 Since the multiplicities of 
edges do not influence of being a loaded tree, we know that $T_1$ and $T_2$ are still loaded trees and 
the weights of vertices $u$ and $v$ stay unchanged after the edge-cutting. 
From \cite[Proposition 5.10]{vertex_splitting} we know that $|\int(T)|=|\int(T_1)|\cdot |\int(T_2)|$.

When $r>1$,
let $s_1$ be the number of edges of $T_1$ and $s_2$ be that of $T_2$.
Let $N:=I_1\cup I_2$ and $|N|=n$. A quick calculation reveals that it can never happen that both
$T_1$ and $T_2$ are proper, which indicates that $\int(T)$ is always zero, if we want an analogous relation
as in the single-edge cutting case.
Therefore, we need to modify our construction.
We construct $T'_1$ (from $T_1$) by removing the label 
$x$ and attaching to $u$ a new vertex $u'$ via an edge $e'$ connecting $u$ and $u'$, where the labeling set 
of $u'$ is $\{a,b\}$ and the multiplicity of $e'$ is set to be $|I_1|+2-s_1-3=|I_1|-s_1-1$; note that $a$, $b$ are 
two new labels not in $N$. The construction
of $T'_2$ via $T_2$ is done analogously. In this way, we are able to obtain proper trees by adding new edges at the 
end of the obtained trees.
Cutting off an edge of $T$, obtaining $T'_1$ and $T'_2$ as stated above, is called a {\bf multi-edge cutting} operation.
In this section, a main task for us is to investigate the relation between the value of $T$ and the values
of $T'_1$ and $T'_2$. In the sequel, we will introduce more algebraic notations, so as to express better the theorem.

First, notice that $T'_1$ and $T'_2$ are both in many cases strictly smaller than $T$, in those cases their monomials live in 
different ambient Chow ring than $T$. For this, we introduce a foot index for the integral symbol,
indicating the ambient space. For instance, let $M$ be a monomial in $A^{n-3}(n)$, then $M$ is in the Chow ring of
$\mathcal{M}_n$. Let $N:=\{1,\cdots,n\}$ be the labeling set of $\mathcal{M}_n$, we denote by 
$\int_{\mathcal{M}_N}(M)$ for the value of $M$. 
In this section, we will use this notation, so as to clarify the ambient Chow ring and the labels in the ambient variety
as well. In order to avoid crashes of notations,
we denote by $A^{\bullet}(X)$ instead of $A^*(X)$ for the Chow ring of the variety $X$, in this section.

Let $M_T$ be the monomial of $T$, let $\delta_{I_1,I_2}$ be the corresponding factor of the edge $e$ and let $r\geq 1$ 
be the multiplicity of 
$e$; note that $I_1\cup I_2=N$. 
Let $\mu_1$ be the product of factors of 
edges in Component-$I_1$ and let $\mu_2$ be that in Component-$I_2$ ---  
note that all factors inherit
the multiplicities of their corresponding edges via their powers. From the property of $M_T$ being a tree 
monomial, we obtain the following conclusion.
Let $S$ be the multi-set of the factors of $M_T$, note that factors with power higher than one appear more than once in the set. Then 
$\mu_1$ is the product of the generators $\delta_{U,V}\in S$ such that $U\subsetneq I_1$, and
$\mu_2$ is the product of the generators $\delta_{U,V}\in S$ such that $V\subsetneq I_2$.
Then $M_T=\mu_1\cdot \mu_2\cdot \delta^r_{I_1,I_2}$ and thence $\int(T)=\int_{\mathcal{M}_N}{\mu_1\cdot \mu_2\cdot \delta^r_{I_1,I_2}}$.
 Let $s_1$, $s_2$ be the degrees of $\mu_1$ and $\mu_2$, respectively; then $\mu_1\in A^{s_1}(\mathcal{M}_N)$, 
 $\mu_2\in A^{s_2}(\mathcal{M}_N)$. Since $T$ is a proper loaded tree, we know that $r=|N|-s_1-s_2-3$.

 Replacing each factor $\delta_{U,V}$ in $\mu_1$ by $\delta_{U, V\setminus I_2\cup \{x\}}$, we obtain the product
 $\gamma_1\in A^{s_1}(\mathcal{M}_{I_1\cup\{x\}})$. 
 Replacing each factor $\delta_{U,V}$ in $\mu_2$ by $\delta_{U, V\setminus I_1\cup \{x\}}$, we obtain the product
 $\gamma_1\in A^{s_2}(\mathcal{M}_{I_2\cup\{x\}})$. 
 With out loss of generality, assume that the labels of $T_1$ form $I_1$ and
 those of $T_2$ form $I_2$. Then it is not hard to see that $M_{T_1}=\gamma_1$ and $M_{T_2}=\gamma_2$. 
 Replacing each factor $\delta_{U,V}$ in $\mu_1$ by $\delta_{U,V\setminus I_2\cup\{a,b\}}$, we get the 
 product $\nu_1\in A^{s_1}(\mathcal{M}_{I_1\cup\{a,b\}})$. 
 Replacing each factor $\delta_{U,V}$ in $\mu_2$ by $\delta_{U,V\setminus I_1\cup\{a,b\}}$, we get the 
 product $\nu_2\in A^{s_2}(\mathcal{M}_{I_2\cup\{a,b\}})$. 
 It is not hard to see that $M_{T'_1}=\nu_1\cdot (\delta_{I_1,\{a,b\}})^{|I_1|-s_1-1}$ and 
 $M_{T'_2}=\nu_2\cdot (\delta_{I_2,\{a,b\}})^{|I_2|-s_2-2}$.
 Then the following theorem will reveal to us the relation between the value of $T$ and the values of
 $T'_1$ and $T'_2$. 
 \begin{theorem}\label{thm:reveal}
 With the notations above, the following equation holds.
  \begin{flalign*}
\int_{\mathcal{M}_N}{\mu_1\cdot \mu_2\cdot (\delta_{I_1,I_2})^r} =  
&{r-1 \choose |I_1|-s_1-2}\cdot \int_{\mathcal{M}_{I_1\cup\{a,b\}}}{\nu_1\cdot (\delta_{I_1,\{a,b\}})^{|I_1|-s_1-1}}&&\\
                      &\cdot \int_{\mathcal{M}_{I_2\cup\{a,b\}}}{\nu_2\cdot (\delta_{I_2,\{a,b\}})^{|I_2|-s_2-1}}&&
\end{flalign*}
 \end{theorem}
 
 \subsection{Correctness proof}
 In this section, we prove the correctness of the forest algorithm given that Theorem~\ref{thm:reveal} holds.
 
 Let $LT=(V,E,h,m)$ be a proper loaded tree with labeling set $N=\{1,\cdots, n\}$; denote by $w$ its weight function. 
First, we cut off all single edges of $LT$; 
 \cite[Lemma 5.12]{vertex_splitting} says that $|\int(LT)|$ is the product of absolute value of 
 all new trees obtained after the series of operations. In the forest algorithm, this operation is equivalent 
 to {\bf deleting all weight-zero vertices of the redundancy tree $RT$ of $LT$ that come from an edge of $LT$}.
 
 Let $l\in V$ be any 
 leaf that has non-zero weight and let $e=\{l,l_1\}$ be its unique incident edge. Now we
 cut off the multiple edge $e$, obtain two new trees $T'_1$ and $T'_2$ and apply the above formula. 
 We see that by definition, in this case, $r-1=m(e)-1=w(e)$. Since $T$ is proper, we have $|I_1|+|h(l)|-3=m(e)+s_1$.
 By definition we have $|h(l)|+\deg(l)-3=w(l)$. Hence we get $|I_1|-s_1-2=m(e)-1-w(l)=w(e)-w(l)$.
 Hence the binomial coefficient on the right hand side of the formula is ${w(e)\choose w(e)-w(l)}={w(e)\choose w(l)}$.
 
 Note that $|I_2|=|h(l)|$, $s_2=0$ and $w(l)=|I_2|+1-3$ (if $T'_2$ is proper), we obtain $|I_2|-s_2-1=w(l)+1$. Hence the new edge $e'_2$
 added to $l$ in $T'_2$ has multiplicity $w(l)+1$, hence its weight is $w(l)$. Hence $T'_2$ is a proper loaded 
 tree with two vertex connecting by an edge and the weights of two vertices are $w(l)$ and $0$, respectively. 
 Since $w(l)\neq 0$, it is a sun-like tree, by \cite[Theorem 7.1]{vertex_splitting}, we know that its absolute
 value is ${w(l)\choose w(l)}=1$. Tree $T'_1$ is obtained by replacing vertex $l$ by a weight-zero vertex, and 
 replacing edge $e$ by an edge $e'_1$ with multiplicity $w(e)-w(l)+1$ and hence the weight of $e'_1$ is $w(e)-w(l)$.
 Hence $\int(T)={w(e)\choose w(l)}\cdot\int(T_1)$. 
 Notice that if some leaf has weight bigger than its parent vertex, after cutting off its unique incident edge, 
 the obtained tree where this leaf lives will be unproper, which leads to the input loaded tree being value zero.
 In the forest algorithm, this step is equivalent to 
 {\bf applying the recursive formula first to all the non-zero-weight leaves of the redundancy forest $RF$ if the leaf comes
 from some leaf of $LT$ and return zero if some leaf in $RF$ has weight bigger than its parent vertex and this leaf
 vertex comes from some leaf of $LT$.} 
 
 We can repeat the above process to all the leaves that have non-zero weight, namely cut off the edges of which the 
 incident leaf has non-zero weight. One can check that the above formula will not help when the leaf has weight zero,
  hence the next step we will do is to consider a {\em star-cut}. 
 
 Note that the concept of star-cut is compatible with multi-edge cut; these two operations can be combined.
Our next step is to do a star-cut, then apply the formula in Theorem \ref{thm:reveal}.
Denote by $T_2$ the star-shaped loaded tree, and by $T_1$ the other loaded tree we obtain from $T$, after
a star-cut. First, it is clear that after the first step, all leaves have weight zero. Hence all leaves of $T_2$ 
have weight zero. Suppose that $T_2$ is proper.
If the middle vertex $v$ of $T_2$ has weight zero, then by its properness, we know that all edges have weight zero as well;
it is a clever tree hence has value $1={0\choose 0,\cdots,0}$. If
its middle vertex has non-zero
weight, then its a sun-like tree, then by \cite[Theorem 7.1]{vertex_splitting} we obtain that 
$|\int(T_2)|={k\choose m_1,\cdots,m_q}$, where $k$ is the weight of $v$ and $m_1,\cdots, m_q$  are the weights of its edges.
One property of multinomial coefficients says that 
$${k\choose m_1,\cdots, m_q}={k\choose m_1}\cdot {k-m_1\choose m_2}\cdot \ldots \cdot {m_q\choose m_q}.$$
Let $e_1$, $e_2$ be the new edges added to $T_1$, $T_2$ respectively, w.l.o.g. let $m_q$ be the weight of 
$e_2$, denote by $w_i$ the weight function of $T_i$, $i=1,2$.
Since $T_2$ is proper, we know that $m_q$ equals the multiplicity of $e_2$ minus $1$, that is, 
$m_q=|I_2|-s_2-2$. Also, we see that $w_1(e_1)=|I_1|-s_1-2$ and that $w_1(e_1)+w_2(e_2)=w(e)$. 
If $T_2$ is not proper, it means that $k\neq \sum_{i=1}^q{m_i}$. Then when we focus on the vertex $w$ in $RF$ that 
comes from this middle vertex, 
recursively apply the forest algorithm, at some point, the weight of $w$ must be less than that of one of its adjacent leaves; or,
when $w$ itself becomes a leaf, its weight is bigger than the parent. Both cases of course should give us zero in the output.
This fact, in the forest algorithm, is described as:
{\bf return zero if some leaf in $RF$ has weight bigger than its parent vertex and this leaf
 vertex does not come from some leaf of $LT$.}

We can repeat the star-cut operation until the loaded tree has only two vertices. 
Actually, any weight-zero vertex will become a middle vertex 
of some star at some stage of the process. From the above analysis,
we see that those zero-weights (no matter on the leaves or in the middle) do not influence or contribute anything to 
our computation, therefore they can be omitted. Based on this fact, the following steps are introduced in the forest algorithm: 
{\bf deleting all weight-zero vertices of the redundancy tree $RT$ of $LT$ that come from a vertex of $LT$
and deleting the leaf from the forest after each recursive step}.

When the loaded tree has only two vertices, we can apply again the multi-edge cut, finally obtain two sun-like trees;
actually the values of the obtained sun-like trees coincides with the base cases.
The calculations of the star and the sun-like trees indicate the recursive formula applied to the redundancy forest.
This concludes the correctness of the forest algorithm.
 
 \subsection{From algebra to geometry}
 In this section, we prove Theorem \ref{thm:reveal} which indicates the main geometric structure hidden beneath the forest algorithm,
 however using pure algebra.

 Let $f:X\to Y$ be a map between two smooth projective varieties. Then $f$ induces the pushforward map
 $f_*:A^{\bullet}(X)\to A^{\bullet}(Y)$, which is a group homomorphism, and the pullback map 
 $f^*:A^{\bullet}(Y)\to A^{\bullet}(X)$, which is a ring homomorphism that preserves the degree of the ambient group where the 
 element lives. Let $\alpha\in A^{\bullet}(X)$, $\beta\in A^{\bullet}(Y)$, then the following adjoint formula on the integrals holds:
 $\int_X(\alpha\cdot f^*(\beta)) = \int_Y(f_*(\alpha)\cdot \beta)$.
 
 It is known that $D_{I_1,I_2}\cong \mathcal{M}_{I_1\cup\{x\}}\times \mathcal{M}_{I_2\cup\{x\}}$, where $x\notin N$ is a new label.
 Denote by $p_{I_1,x,N}$ the projection from $D_{I_1,I_2}$ to $\mathcal{M}_{I_1\cup\{x\}}$
 and by $p_{I_2,x,N}$ the projection from $D_{I_1,I_2}$ to $\mathcal{M}_{I_2\cup\{x\}}$.
Denote by $i_{I_1,I_2}$ the embedding of $D_{I_1,I_2}$ as a hypersurface into $\mathcal{M}_N$.
 Let $i:=i_{I_1,I_2}$ and let $p_1:=p_{I_1,x,N}$, $p_2:=p_{I_2,y,N}$.
 Then we have $i^*(\mu_1\cdot \mu_2)=p^*_1(\gamma_1)\cdot p^*_2(\gamma_2)$. 
 Then apply the pushforward map on both sides, we obtain $\mu_1\cdot \mu_2\cdot\delta_{I_1,I_2}=i_*(p_1^*(\gamma_1)\cdot p_2^*(\gamma_2))$;
 this equation will be used later in the proof of Theorem~\ref{thm:reveal}.
Define $\beta_{x,I_1\cup\{x\}}:=[(p^{-1}_{I_1,x,I_1\cup \{a,b\}})^*\circ i^*_{I_1,\{a,b\}}](\delta_{I_1,\{a,b\}})$.
Note that $p_{I_1,x,I_1\cup\{a,b\}}$ is the projection from $D_{I_1,\{a,b\}}$ to $\mathcal{M}_{I_1\cup \{x\}}$.
Consider the isomorphism $D_{I_1,\{a,b\}}\cong \mathcal{M}_{I_1\cup\{x\}}\times \mathcal{M}_{\{a,b,x\}}$; 
since $\mathcal{M}_{\{a,b,x\}}$ is just a point, we obtain that $D_{I_1,\{a,b\}}\cong \mathcal{M}_{I_1\cup \{x\}}$.
Therefore the inverse of $p_{I_1,x,I_1\cup\{a,b\}}$ exists. Since pullback is degree-preserving and 
$\delta_{I_1,\{a,b\}}\in A^1(\mathcal{M}_{I_1\cup\{a,b\}})$, we know that $\beta_{x,I_1\cup\{x\}}\in A^1(\mathcal{M}_{I_1\cup\{x\}})$. 
Analogously, we define $\beta_{x,I_2\cup\{x\}}\in A^1(\mathcal{M}_{I_2\cup\{x\}})$, simply by 
replacing $I_1$ by $I_2$, in the definition of $\beta_{x,I_1\cup\{x\}}$. 
Let $\beta_1:=\beta_{x,I_1\cup\{x\}}$, $\beta_2:=\beta_{x,I_2\cup\{x\}}$.
\begin{lemma}\label{lem:direct_product}
  The following equation holds:
 $$i^*(\delta_{I_1,I_2})=p^*_1(\beta_1)+p^*_2(\beta_2),$$
 where $x\notin N$ is a new label.
\end{lemma}
In order to prove the above lemma, we need to introduce some basic properties of the Chow group of 
a direct product of two varieties in $\mathcal{M}_N$. 
Let $X$ and $Y$ be two smooth projective varieties of $\mathcal{M}_N$, then we have 
$A^1(X\times Y)\cong A^1(X)\bigoplus A^1(Y)$.
Let $\pi_l$, $\pi_r$ be the projection from 
$X\times Y$ to $X$ and $Y$, respectively. We know that for any $y_0\in Y$, there exists a right inverse 
$\sigma_l$ of $\pi_l$ such that $\sigma_l(x):=(x,y_0)$. Let $\sigma_l$ be any such inverse; the choice 
of the element in $Y$ does not matter; we define $\sigma_r$ analogously, as a right inverse for $\pi_r$.
Then for any $t\in A^1(X\times Y)$, we have $t=\pi^*_l\circ \sigma_l^*(t)+\pi^*_r\circ \sigma_r^*(t)$.
Observe that we have the isomorphism $D_{I_1,I_2}\cong \mathcal{M}_{I_1\cup\{x\}}\times \mathcal{M}_{I_2\cup\{x\}}$
in $\mathcal{M}_N$. Let $q_{I_i,x,N}$ be any right inverse (as described above) of $p_{I_i, x,N}$ for $i=1,2$. 
Denote by $q_i:=q_{I_i,x,N}$, for $i=1,2$.
Then, from the above analysis, we know that for any $t\in A^1(D_{I_1,I_2})$, we have 
$t=p_1^*\circ q_1^*(t)+p_2^*\circ q_2^*(t)$. Now we will introduce an equation telling us the 
relation between the integral value on the direct product and the integral values on its two coordinates;
this will be used later in our proof.
Let $\alpha\in A^{\bullet}(X)$, $\beta\in A^{\bullet}(Y)$.
Then we have 
$$\int_{X\times Y}{\pi_l^*(\alpha)\cdot \pi_r^*(\beta)} = \int_X{\alpha}\cdot \int_Y{\beta}.$$

We need one more preparation before proving Lemma \ref{lem:direct_product}.
Define $s_{k,l,N}:\mathcal{M}_{N\setminus \{k\}}\to \mathcal{M}_N$
as $$s_{k,l,N}:=i_{\{k,l\},N\setminus \{k,l\}}\circ p^{-1}_{N\setminus\{k,l\},l,N},$$ 
where $k,l$ are two distinct labels of $N$. 
Note that $p_{N\setminus\{k,l\},l,N}$ is an isomorphism, hence it has an inverse. 
There is a surjective forgetful map $c_{a,N}:\mathcal{M}_N\to \mathcal{M}_{N\setminus\{a\}}$
for any $a\in N$. The above defined map $s_{k,l,N}$ is a right inverse of $c_{k,N}$. The image 
of $s_{k,l,N}$ is the hypersurface $D_{\{k,l\},N\setminus \{k,l\}}$ in $\mathcal{M}_N$.

\begin{proof}[Proof of Lemma \ref{lem:direct_product}]
Recall that $i^*$ is the pullback map from $A^{\bullet}(\mathcal{M}_N)$
 to $A^{\bullet}(D_{I_1,I_2})$ and that $\delta_{I_1,I_2}\in A^1(\mathcal{M}_N)$.
 Since pullback is a degree-preserving ring homomorphism, we know that
 $i^*(\delta_{I_1,I_2})\in A^1(D_{I_1,I_2})$. Using the result from earlier
 analysis, we have Equation (a):
 $$i^*(\delta_{I_1,I_2})=p^*_1\circ q_1^*(i^*(\delta_{I_1,I_2}))+ p^*_2\circ q_2^*(i^*(\delta_{I_1,I_2})).$$
 We claim that it suffices to prove 
 Equation (b): $q^*_1\circ i^*(\delta_{I_1,I_2})=\beta_1$ 
 and Equation (c): $q^*_2\circ i^*(\delta_{I_1,I_2})=\beta_2$. Suppose they hold,
 then from (b) we have: $p_1^*\circ q_1^*\circ i^*(\delta_{I_1,I_2})=p_1^*(\beta_1)$. 
 Analogously, we obtain $(q_2\circ p_2)^*(i^*(\delta_{I_1,I_2}))=p^*_2(\beta_2)$
 from (c). Substituting the equalities back to Equation (a), we obtain the wanted 
 equation. Since (b) and (c) are symmetric, it suffices to prove (b).
 
 We prove Equation (b) by induction on $|I_2|$. Recall the definition of $\beta_1$, in this case, we have:
 $\beta_1=(p_1^{-1})^*\circ i^*(\delta_{I_1,I_2})$. It suffices to show that $(p_1^{-1})^*=q_1^*$. 
 Since $q_1$ is a right inverse of $p_1$, we have $p_1\circ q_1=id$.
 But in this case $p_1$ is an isomorphism, so does $q_1$. Therefore, $q_1=p_1^{-1}$. Hence $(p_1^{-1})^*=q_1^*$.
 Suppose Equation (b) holds when $|I_2|=z-1$, now $|I_2|=z\geq 3$, $z\in \mathbb{N}$.
 Let $k,l\in I_2$ be two distinct labels and let $I'_2:=I_2\setminus \{k\}$, $N':=N\setminus \{k\}=I_1\cup I'_2$.
 Then we have the following equality: $i\circ q_1=s_{k,l,N}\circ i_{I_1,I'_2}\circ q_{I_1,x,N'}$ of maps from 
 $\mathcal{M}_{I_1\cup \{x\}}$ to $\mathcal{M}_N$. Since pullback is a contravariant functor and
 $\delta_{I_1,I_2}\in A^{\bullet}(\mathcal{M}_N)$, we obtain:
 $q^*_1\circ i^*(\delta_{I_1,I_2})=q^*_{I_1,x,N'}\circ i^*_{I_1,I'_2}\circ s^*_{k,l,N}(\delta_{I_1,I_2})
 =q^*_{I_1,x,N'}\circ i^*_{I_1,I'_2}(\delta_{I_1,I'_2})$. Now we can use the induction hypothesis,
 since $|I'_2|=|I_2|-1=z-1$. Hence we have $q^*_1\circ i^*(\delta_{I_1,I_2})=\beta_1$.
\end{proof}

\begin{proof}[Proof of Theorem \ref{thm:reveal}]
 From earlier analysis, we have $$\mu_1\cdot \mu_2\cdot\delta_{I_1,I_2}=i_*(p_1^*(\gamma_1)\cdot p_2^*(\gamma_2)).$$
 Use the adjoint formula between pullback and pushforward, the result in Lemma \ref{lem:direct_product}, 
 the property that the pullback being a ring homomorphism and the integral map being a group homomorphism. Then consider the isomorphism 
 $D_{I_1,I_2}\cong \mathcal{M}_{I_1\cup\{x\}}\times \mathcal{M}_{I_2\cup\{x\}}$, and use the 
fact on direct product of varieties mentioned earlier, we further get
   \begin{flalign*}
\int_{\mathcal{M}_N}{\mu_1\cdot \mu_2\cdot (\delta_{I_1,I_2})^r} &=  
\int_{\mathcal{M}_N}{(\delta_{I_1,I_2})^{r-1}\cdot i_*(p_1^*(\gamma_1)\cdot p_2^*(\gamma_2)) }&&\\
                      &=\int_{D_{I_1,I_2}}{i^*((\delta_{I_1,I_2})^{r-1})\cdot p_1^*(\gamma_1)\cdot p_2^*(\gamma_2)}&&\\
                       &= \int_{D_{I_1,I_2}}{(p^*_1(\beta_1)+p^*_2(\beta_2))^{r-1}\cdot p_1^*(\gamma_1)\cdot p_2^*(\gamma_2)}&&\\
     = \sum_{k=0}^{r-1}&{r-1\choose k}\cdot {\int_{D_{I_1,I_2}}{p_1^*(\beta_1^k\cdot \gamma_1)\cdot p_2^*(\beta_2^{r-1-k}\cdot \gamma_2)}}&&\\
                      = \sum_{k=0}^{r-1}&{r-1\choose k}\cdot \int_{\mathcal{M}_{I_1\cup\{x\}}}{\beta_1^k\cdot \gamma_1} \cdot 
                     \int_{\mathcal{M}_{I_2\cup\{x\}}}{\beta_2^{r-1-k}\cdot \gamma_2}&&
\end{flalign*}
Recall that the integral value is defined to be zero if the monomial is not in the Chow group of codimension $3$ of the ambient Chow ring.
Therefore, we can already omit those summands that are zero in the above sum. Since $\beta_1\in A^1(\mathcal{M}_{I_1\cup\{x\}})$
and $\beta_2\in A^1(\mathcal{M}_{I_2\cup\{x\}})$, we see that we only need to consider the summands such that
$k+s_1=|I_1|+1-3$ and $r-1-k+s_2=|I_2|+1-3$ hold, that is, $k=|I_1|-s_1-2=r+1+s_2-|I_2|$. 
Hence there is only one summand left, we obtain the following
formula:
   \begin{flalign*}
&\int_{\mathcal{M}_N}{\mu_1\cdot \mu_2\cdot (\delta_{I_1,I_2})^r}&&\\
= {r-1\choose |I_1|-s_1-2}&\cdot
\int_{\mathcal{M}_{I_1\cup \{x\}}}{\beta_1^{|I_1|-s_1-2}\cdot \gamma_1}
\cdot \int_{\mathcal{M}_{I_2\cup \{x\}}}{\beta_2^{|I_2|-s_2-2}\cdot \gamma_2}.&&
\end{flalign*}
As a special case, let $I_2=\{a,b\}$, then we have $\mu_1=\nu_1$, $\mu_2=1$, $s_2=0$.
In this case, we see that $|I_1|+2-3=s_1+r$, that is, $r=|I_1|-s_1-1$. The formula becomes
$$\int_{\mathcal{M}_{I_1\cup\{a,b\}}}{\nu_1\cdot (\delta_{I_1,\{a,b\}})^{|I_1|-s_1-1}}
=\int_{\mathcal{M}_{I_1\cup\{x\}}}{\beta_1^{|I_1|-s_1-2}}\cdot \gamma_1,$$
analogously, when $I_1=\{a,b\}$, we get the following equation:
$$\int_{\mathcal{M}_{I_2\cup\{a,b\}}}{\nu_2\cdot (\delta_{I_2,\{a,b\}})^{|I_2|-s_2-1}}
=\int_{\mathcal{M}_{I_2\cup\{x\}}}{\beta_2^{|I_2|-s_2-2}}\cdot \gamma_2.$$
The statement then follows.
\end{proof}
\section{Acknowledgement}
The research was funded by the Austrian Science Fund (FWF): W1214-N15, project DK9.

I am very much grateful to Josef Schicho for formulating my rough graphical ideas down as in Theorem \ref{thm:reveal} (exquisite!)
and providing me with the proof (technical but astounding!) of it, which to a great extend realized the correctness proof, hence further 
drew a perfect full stop on the forest algorithm. I thank Nicolas Allen Smoot for 
helping me express in a better way the proof of the existence of star-cuts.



\begin{thebibliography}{99}
\bibitem{GW}
Behrend Kai.
\newblock Gromov-Witten invariants in algebraic geometry.
\newblock {\em  Inventiones mathematicae}, 127.3 (1997): 601-617.

\bibitem{km1}
Knudsen Finn, and Mumford David.
\newblock The Projectivity of the moduli space of stable curves I: Preliminaries on ``det'' and ``Div''.
\newblock {\em Mathematica Scandinavica}, 39.1 (1977), 19-55.

\bibitem{km2}
Knudsen Finn F.
\newblock The projectivity of the moduli space of stable curves, II: The stacks $M_{g,n}$. 
\newblock {\em Mathematica Scandinavica}, 52.2 (1983), 161-199.

\bibitem{km3}
Knudsen Finn F.
\newblock The Projectivity of the Moduli Space of Stable Curves, III: 
The Line Bundles on $M_{g,n}$, and a proof of the Projectivity of $M_{g,n}$ in Characteristic $0$. 
\newblock {\em Mathematica Scandinavica} (1983), 200-212.


\bibitem{keel}
Keel Sean.
\newblock Intersection theory of moduli space of stable $n$-pointed curves of 
genus zero.
\newblock {\em Transaction of the American Mathematical Society}, {\bf 330} (1992),
no. 2, 545-574.


\bibitem{laman}
Matteo Gallet, Georg Grassegger, Josef Schicho.
\newblock Counting realizations of Laman graphs on the sphere.
\newblock {\em The Electronic Journal of Combinatorics}, Volume 27, Issue 2 (2020).

\bibitem{QS:representation}
Qi Jiayue, Schicho Josef.
\newblock Five Equivalent Ways to Describe Phylogenetic Trees.
\newblock {\em arXiv:2011.11774}, preprint.

\bibitem{ACM_communication}
Qi Jiayue.
\newblock A calculus for monomials in Chow group of zero cycles in the moduli 
space of stable curves.
\newblock {\em DK-Report}, No. 2020-11, September 2020.

\bibitem{vertex_splitting}
Qi Jiayue.
\newblock A tree-based algorithm on monomials in the Chow group of zero cycles
in the moduli space of stable pointed curves of genus zero.
\newblock {\em arXiv:2101.03789}, preprint.

\end{thebibliography}
\end{document}